\documentclass[final,leqno,onefignum,onetabnum]{siamltex}

\usepackage{subfigure,upgreek,dsfont,braket,amssymb,amsmath}
\usepackage[small,bf]{caption}
\usepackage{booktabs,dsfont}
\usepackage{algorithm}
\usepackage{float}
\usepackage{graphicx}
\usepackage{algorithmicx}
\usepackage{color,float}
\usepackage{url}

 \newtheorem{remark}[theorem]{Remark}

\newcommand{\de}{\delta}

\newcommand{\auskommentieren}[1]{}

\newcommand{\beq}{\begin{equation}}
\newcommand{\eeq}{\end{equation}}

\title{Alternative to evolving surface finite element method} 

\author{Maryia Borukhava \thanks{Bronnenmayerstraße 20, 73037 G\"oppingen, Germany, {\tt  Maryia.Borukhava@studium.uni-hamburg.de}}
\and Heiko Kr\"oner\thanks{Schwerpunkt Optimierung und Approximation, Universit\"at Hamburg, Bundesstra\ss e 55, 20146 Hamburg, Germany.
{\tt  heiko.kroener@uni-hamburg.de}}}

\begin{document}
\maketitle
\slugger{mms}{xxxx}{xx}{x}{x--x}

\begin{abstract}
ESFEM is a method introduced in \cite{DE07} in order to solve a linear advection-diffusion equation on an evolving two-dimensional surface with finite elements by using a moving grid with nodes sitting on and evolving with the surface.  The evolution of the surface is assumed to be given as a smooth one-parameter family of embeddings of a fixed initial surface into $\mathbb{R}^3$ satisfying uniform $C^4$ bounds. 
We calculate an equivalent transformed equation which is defined on the fixed initial surface and can hence be solved numerically on a fixed grid. 
We present numerical examples which indicate that both approaches are essentially of the same accuracy.
\end{abstract}


\pagestyle{myheadings}
\thispagestyle{plain}
\markboth{A remark on the evolving surface finite element method}{Maryia Borukhava and Heiko Kr\"oner}

\section{Introduction} \label{intro}
In many applications it is important to consider PDEs which are defined on surfaces and not in Euclidean space, especially in the case of parabolic equations it is of interest to assume that these surfaces (where the equation is defined) evolve with respect to time in a certain prescribed way.
In \cite{DE07} the so-called evolving surface finite element method (ESFEM) is proposed in order to solve an advection-diffusion equation on an evolving surface, cf. \cite[Sections 1.1 and 1.2]{DE07}.
This setting models e.g. the transport of an insoluble surfactant on the interface between two flowing fluids or pattern formation on the surfaces of growing organisms modeled by reaction-diffusion equations, cf. \cite[Section 1.4]{DE07} for further and a more detailed exposition of applications. 
There are several papers which deal with linear parabolic equations on evolving surfaces, e.g. in \cite{KG, DE13} it is shown that classical $L^2$- and $L^{\infty}$- estimates 
carry over to ESFEM and in \cite{LMV13, DLM12} it is shown that this also holds for error 
estimates of Runge-Kutta schemes and Backward difference schemes; we also mention 
\cite{OR2014, ORX2014}.

We present the idea of ESFEM according to \cite{DE07} and note that this method is
introduced therein in order to solve a special linear parabolic equation, namely Equation
(\ref{2}), with the finite element method.

Let $\Gamma_0$ be a smooth, closed, connected and oriented hypersurface in $\mathbb{R}^3$. 
Let
$\Phi(t, \cdot): \Gamma_0 \rightarrow \mathbb{R}^3$, $t \in [0, T_0]$, $T_0>0$, be a
family of embeddings, $\Phi$ smooth, $\Gamma(t) = \Phi(t,\Gamma_0)$ the moving surfaces. We define 
the set 
\beq \label{1}
G_{T_0} = \bigcup_{t \in [0, T_0]} \{t\}\times \Gamma(t)
\eeq
and consider there the advection-diffusion equation
\beq \label{2}
\dot u + u \nabla^{\Gamma(t)} \cdot v - \nabla^{\Gamma(t)}\cdot (D_0 \nabla^{\Gamma(t)}u) =0
\eeq
where $\dot u$ denotes the material derivative, $v(p_t, t)= \frac{\partial \Phi}{\partial t}(\Phi(t,\cdot)^{-1}(p_t), t)$ 
the velocity of the moving surface in a point $(t,p_t) \in G_{T_0}$ and hence 
$\nabla^{\Gamma(t)} \cdot v $ its tangential divergence. Furthermore, the diffusion coefficient $D_0: G_{T_0}\rightarrow 
Mat(n+1, \mathbb{R})$ is so that it vanishes on the normal space of the moving surfaces,
i.e. $D_0(t,p_t) \nu=0$ for all $\nu \in N_{p_t}$ where $N_{p_t}$ is the normal space of $\Gamma(t)$ 
in $p_t$.
Equation (\ref{2}) can be written in variational form as
\begin{equation} \label{w10}
\frac{d}{dt}\int_{\Gamma(t)}u\varphi + \int_{\Gamma(t)}D_0 \nabla^{\Gamma(t)}\cdot \nabla^{\Gamma(t)}\varphi =
\int_{\Gamma(t)}u \dot \varphi \quad \forall \varphi \in C^{\infty}(G_{T_0}).
\end{equation}
ESFEM  considers at every time $t$ an interpolating polyhedral surface $\Gamma_h(t)$ which approximates the evolving surface $\Gamma(t)$ and which consists of triangles with vertices $X_j(t)=\Phi(t,X_j(0))$, $j=1, ..., N$, sitting on $\Gamma(t)$ and moving with the surface. Here, $X_j(0)$, $j=1, ..., N$, are fixed nodes on the initial surface $\Gamma(0)$. The finite element basis functions $\varphi_i$, $i=1, ..., N$,  are
defined on
\beq \label{1__}
G^h_{T_0} = \bigcup_{t \in [0, T_0]} \{t\} \times \Gamma_h(t)
\eeq
and chosen so that $\varphi_i(t, \cdot)$ is piecewise linear (i.e. linear on each triangle of $\Gamma_h(t)$) with 
$\varphi_i(t, X_j(t))= \de_{ij}$. 

 The fully discrete scheme from \cite{DE07} in order to solve (\ref{2}) which uses 
 ESFEM can be found in \cite[Equation (7.2)]{DE07} and is as follows.
 Let $t_0=0< ... < t_M=T_0$, $M \in \mathbb{N}$, be a partition of the time interval, 
 $\Gamma_h^m = \Gamma_h(t_m)$, $u^m=u(\cdot, t_m)$ for a function $u$ on $G^h_{T_0}$ 
 and $u_h^0\in \Gamma_h(0)$ an initial function (and for simplicity $D_0$ the identity) 
 then we solve for $m=0, ..., M-1$ the linear system
 \begin{equation}
 \begin{aligned}
 \frac{1}{\tau}&\int_{\Gamma_h^{m+1}}u_h^{m+1}\varphi^{m+1}_j + \int_{\Gamma_h^{m+1}}\nabla^{\Gamma_h^{m+1}}
 u_u^{m+1}\cdot \nabla^{\Gamma_h^{m+1}}\varphi_j^{m+1} \\
 =& \frac{1}{\tau}\int_{\Gamma_h^m}u_h^m\varphi_j^m, \quad j=1, ..., N.
 \end{aligned}
 \end{equation}

Our different approach to solve (\ref{2}) with the finite element method 
reformulates Equation (\ref{2}) as an equivalent equation on the  
surface $\Gamma(0)$, cf. Equation (\ref{n10}), and uses a fixed grid (the one consisting
of the nodes $X_j(0)$, $j=1, ..., N$,) for the finite element approximation 
(again an implicit Euler method) which leads to a finite element solution 
$\tilde u_h^m$ defined on $\Gamma_h(0)$. Obviously, both approaches coincide if 
$\Phi(t, \cdot)$ is the identity, i.e. if there is no motion.

Our aim is to use both approaches, i.e ESFEM and the method which uses the transformed equation,  in some example cases  and to provide the corresponding error tables. 

The paper is organized as follows. In Section \ref{sec1} we recall some facts about hypersurfaces in the Euclidean space, in Section \ref{reform_fixed_mesh} we derive the reformulated equation, in Section \ref{Numerical_Examples} we present our chosen examples and a more detailed description of the parts from which the coefficients of the transformed equation are put together in the implementation. The error tables of the numerical calculations can be found in Tables 1 to 6.
 
\section{Hypersurfaces in $\mathbb{R}^{n+1}$} \label{sec1}
We recall some facts and notations of embedded hypersurfaces in $\mathbb{R}^{n+1}$ from \cite{E04}. 
 Let $F: \Omega \rightarrow \mathbb{R}^{n+1}$
with $\Omega \subset \mathbb{R}^n$ open be a smooth embedding and $M=F(\Omega)$. For $p \in \Omega$ the coordinate tangent vectors
$\partial_iF(p) = \frac{\partial F}{\partial p_i}(p)$, $1 \le i \le n$, provide a basis of the tangent
space $T_xM$ at $x=F(p)$. 

The metric on $M$ is given by
\beq
g_{ij} = \partial_iF \cdot \partial_jF
\eeq
for $1 \le i,j \le n$, the inverse metric by $(g^{ij})=(g_{ij})^{-1}$.

The tangential gradient of a function $h: M \rightarrow \mathbb{R}$ is defined by
\beq
\nabla^Mh = g^{ij}\partial_jh \partial_iF
\eeq
where we sum over repeated indices.

For a smooth tangent vector field $X = X^i\partial_i F=g^{ij}X_j\partial_iF$ on $M$ (note that 
$X_i = X \cdot \partial_iF$) we define the covariant derivative tensor by
\beq
\nabla^M_iX^j = \partial_iX^j+ \Gamma_{ik}^jX^k
\eeq
where the Christoffel symbols $\Gamma_{ij}^k$ are given by
\beq
\Gamma_{ij}^k = \frac{1}{2}g^{kl}(\partial_ig_{jl}+\partial_j g_{il}-\partial_l g_{ij}).
\eeq
The tangential divergence of $X$ on $M$ is defined by
\beq
\nabla^M\cdot X = \operatorname{div}_MX = \nabla^M_iX^i
\eeq
and the Laplace-Beltrami operator of $h$ on $M$ by
\beq
\Delta_Mh = \operatorname{div}_M\nabla^Mh = \nabla^M\cdot(\nabla^Mh).
\eeq

For a smooth vector field $X: M \rightarrow \mathbb{R}^{n+1}$ which is not necessarily tangent on $M$ we can
also define the divergence with respect to $M$ by
\beq \label{8}
\nabla^M \cdot X = \operatorname{div}_MX = g^{ij}\partial_iX\cdot \partial_jF
\eeq
which reduces to the above expression if $X$ is tangent on $M$. We remark that, of course,  
the divergence in (\ref{8}) is (and
transforms like) a scalar function (when changing local coordinates of $M$). Furthermore, the 
tangential gradient transforms like a scalar function when changing coordinates in $M$.

\section{Reformulation of the equation on a fixed surface} \label{reform_fixed_mesh}
In this section we derive Equation (\ref{n10}) which is an equivalent reformulation of (\ref{2}) on $\Gamma(0)$. The calculation is straightforward but we still present details.

Let $\Omega, \tilde \Omega \subset \mathbb{R}^n$ be open and $\Phi: \tilde \Omega \rightarrow \Omega$
a diffeomorphism. The linear differential operator $L:H^2(\Omega)\rightarrow L^2(\Omega)$
\begin{equation}
Lu = a^{ij}D_iD_j u+ b^i D_i u + c u 
\end{equation}
with coefficients $a^{ij}, b^i, c \in L^{\infty}(\Omega)$, $a^{ij}$ symmetric, transforms into the operator 
$\tilde L:H^2(\tilde \Omega)\rightarrow L^2(\tilde \Omega)$
\begin{equation}
\tilde L \tilde u = \tilde a^{ij}D_i D_j \tilde u + \tilde b^i D_i \tilde u + \tilde c \tilde u, 
\end{equation}
i.e. we have $(Lu) \circ \Phi = \tilde L \tilde u$
for the quantity $\tilde u=u\circ \Phi$. Here, by writing $x(\tilde x) = \Phi(\tilde x)$ and $\tilde x(x) = \Phi^{-1}(x)$
we set 
\begin{equation} \label{n1}
\begin{aligned}
\tilde a^{ij} =& a^{rs} \frac{\partial \tilde x^i}{\partial x^r}\frac{\partial \tilde  x^j}{\partial x^s} \\
\tilde b^k =& b^m \frac{\partial \tilde x^k}{\partial x^m}-a^{ij}\frac{\partial^2  x^m}{\partial \tilde x^r \partial \tilde x^s}
\frac{\partial \tilde x^r}{\partial  x^i}\frac{\partial \tilde x^s}{\partial x^j}\frac{\partial \tilde x^k}{\partial x^m} \\
\tilde c =& c 
\end{aligned}
\end{equation}
and $\circ \Phi$ operations are suppressed. 

These formulas carry over to the surface case when using local coordinates. Let now $\Omega, \tilde \Omega$ be (open subsets of) hypersurfaces in $\mathbb{R}^3$, $\Phi: \tilde \Omega \rightarrow \Omega$ a diffeomorphism and let $a^{ij}$, $b^i$ and $c$ be $L^{\infty}$-sections of the tensor bundles $T^{2,0}(\Omega)$, $T^{1,0}(\Omega)$ and $T^{0,0}(\Omega)$, respectively, and assume that $a^{ij}$ is symmetric. Let  $L$ be defined by
\begin{equation}
Lu = a^{ij}\nabla^{\Omega}_i \nabla^{\Omega}_j u+ b^i \nabla^{\Omega}_i u + c u, \quad u \in H^2(\Omega),
\end{equation}
where $\nabla^{\Omega}$ denotes the Levi-Cevita connection with Christoffel symbols $\Gamma^k_{ij}$ on $\Omega$  (and $\nabla^{\tilde \Omega}$ and $\tilde \Gamma^k_{ij}$ correspondingly on $\tilde \Omega$) then we have
\begin{equation}
\begin{aligned}
Lu =& a^{ij}D_iD_j u+ (b^k-a^{ij}\Gamma^k_{ij}) D_k u + c u \\
=& a^{ij}D_iD_j u+ \bar b^k D_k u + c u
\end{aligned}
\end{equation}
 where $D_i$ denote ordinary partial derivatives. 
Using the formulas (\ref{n1}) we get a transformed operator 
\begin{equation}
\begin{aligned}
\tilde L \tilde u = \tilde a^{ij}\nabla^{\tilde \Omega}_i \nabla^{\tilde \Omega}_j \tilde u+ \tilde b^i \nabla^{\tilde \Omega}_i \tilde u + \tilde c \tilde u, \quad \tilde u \in H^2(\tilde \Omega),
\end{aligned}
\end{equation}
where now in local coordinates $(x^i)$ of $\Omega$ and $(\tilde x^i)$ of $\tilde \Omega$ we have
\begin{equation} \label{n2}
\begin{aligned}
\tilde a^{ij} =& a^{rs} \frac{\partial \tilde x^i}{\partial  x^r}\frac{\partial \tilde x^j}{\partial x^s} \\
\tilde b^k =& \bar b^m \frac{\partial \tilde x^k}{\partial x^m}-a^{ij}\frac{\partial^2  x^m}{\partial \tilde x^r \partial \tilde x^s}
\frac{\partial \tilde x^r}{\partial x^i}\frac{\partial \tilde x^s}{\partial x^j}\frac{\partial \tilde x^k}{\partial x^m} \\
&+ a^{rs}\frac{\partial  \tilde x^i}{\partial x^r}\frac{\partial \tilde x^j}{\partial x^s} \tilde \Gamma^k_{ij}\\
\tilde c =& c. 
\end{aligned}
\end{equation}
A choice of local coordinates in $\tilde \Omega$ induces via $\Phi$ local coordinates in $ \Omega$ and we stipulate that in the following the local coordinates of $\Omega$ and $\tilde \Omega$ are related in this way. Then the formulas (\ref{n2}) simplify to
\begin{equation} \label{n3}
\begin{aligned}
\tilde a^{ij} =& a^{ij}  \\
\tilde b^k =& \bar b^k 
+ a^{ij}\tilde \Gamma^k_{ij}\\
\tilde c =& c. 
\end{aligned}
\end{equation}
Let us consider the case where the main part is in divergence form
\begin{equation}
\begin{aligned}
Lu =& \nabla^{\Omega}_i(a^{ij} \nabla^{\Omega}_j u)+ b^i \nabla^{\Omega}_i u + c u \\
=& a^{ij}\nabla^{\Omega}_i \nabla^{\Omega}_j u + (\nabla^{\Omega}_ja^{ij}+b^i)\nabla^{\Omega}_i u+cu
\end{aligned}
\end{equation}
then we get the transformed operator with main part in divergence form
\begin{equation} \label{m1}
\begin{aligned}
\hat L \tilde u =& \hat a^{ij}\nabla^{\tilde \Omega}_i \nabla^{\tilde \Omega}_j \tilde u+ \hat b^i \nabla^{\tilde \Omega}_i \tilde u + \hat c \tilde u \\
=& \nabla^{\tilde \Omega}_i (\hat a^{ij}\nabla^{\tilde \Omega}_j \tilde u)+(\hat b^i-\nabla^{\tilde \Omega}_j\hat a^{ij}) \nabla^{\tilde \Omega}_i \tilde u + \hat c \tilde u \\
=& \nabla^{\tilde \Omega}_i (\hat a^{ij}\nabla^{\tilde \Omega}_j \tilde u)+\check b^i  \nabla^{\tilde \Omega}_i \tilde u + \hat c \tilde u
\end{aligned}
\end{equation}
where 
\begin{equation} \label{n4}
\begin{aligned}
\hat a^{ij} =& a^{ij}  \\
\hat b^k =& \bar b^k 
+ a^{ij}\tilde \Gamma^k_{ij}+\nabla^{\Omega}_ia^{ik} =   b^k 
+ a^{ij}(\tilde \Gamma^k_{ij}-\Gamma^k_{ij})+\nabla^{\Omega}_ia^{ik}\\
\hat c =& c
\end{aligned}
\end{equation}
and hence
\begin{equation} \label{n5}
\check b^k =   b^k 
+ a^{ij}(\tilde \Gamma^k_{ij}-\Gamma^k_{ij})+\nabla^{\Omega}_ia^{ik}-\nabla^{\tilde \Omega}_i a^{ik}.
\end{equation}

Let us assume that $u$ and $\Phi$ are as in Section  \ref{intro} and that $D_0$ is the identity.
We define $\tilde u(t,x) = u(t, \Phi(t,x))$ and transform (\ref{2}) in the equivalent equation
\begin{equation} \label{n10}
\begin{aligned}
\frac{d}{dt}\tilde u & -\nabla^{\Gamma(0)}_i(g^{ij}(t)\nabla^{\Gamma(0)}_j \tilde u)  \\
&+ (g^{ij}(t)( \Gamma^k_{ij}(t)-\Gamma^k_{ij}(0))+ \nabla^{\Gamma(0)}_jg^{kj}(t))\nabla^{\Gamma(0)}_k \tilde u
+ \tilde u \nabla^{\Gamma(t)} \cdot v= 0
\end{aligned}
\end{equation}
on $[0, T_0]\times \Gamma(0)$ where the time derivative is now a usual partial derivative, $g_{ij}(t)$ denotes the metric and $\Gamma^k_{ij}(t)$  the 
Christoffel symbols of $\Gamma(t)$. Here, we used that the covariant derivative of the metric vanishes and the coupling of local coordinates via $\Phi(t, \cdot)$.

\begin{remark} \rm
 We note that if one considers instead of (\ref{2}) a general linear parabolic equation
 \begin{equation}
  \dot u -a^{ij}\nabla^{\Gamma(t)}_i \nabla^{\Gamma(t)}_j u + b^i\nabla^{\Gamma(t)}_iu + c u = f
 \end{equation}
 on $G_{T_0}$ where
 \begin{equation}
  a^{ij}: G_{T_0}\rightarrow T^{2,0}(\Gamma(t)), \quad b^i:G_{T_0}\rightarrow T^{1,0}(\Gamma(t)),
  \quad c,f:G_{T_0}\rightarrow \mathbb{R}
 \end{equation}
 so that $a^{ij}(t, \cdot)$ and $b^i(t, \cdot)$ are sections of the corresponding bundles
 then the transformation rules (\ref{n4}) and (\ref{n5}) -- of course --
 provide a reformulation of this equation on the fixed surface as well.
\end{remark}

\section{Examples and implementation} \label{Numerical_Examples}
We choose $\Gamma(0)=\partial B_1(0)\subset \mathbb{R}^3$, $T_0=\frac{3}{2}$ and define the motion of the surface by
\beq
\Phi(t,x) = A(t)x
\eeq
where we consider the cases
\begin{equation} \label{case1}
A(t) = 
\left(
\begin{matrix}
\sqrt{1+5 \sin t} & 0 & 0 \\ 0 & 1 & 0 \\ 0 & 0 & 1
\end{matrix}
\right),
\end{equation}
\begin{equation} \label{case2}
A(t) = 
\left(
\begin{matrix}
1+tx_1^2 & 0 & 0 \\ 0 & 1 & 0 \\ 0 & 0 & 1
\end{matrix}
\right)
\end{equation}
and
\begin{equation} \label{case3}
A(t) = 
\left(
\begin{matrix}
\cos(\eta \mu(x_3)t) & -\sin(\eta \mu(x_3)t) & 0 \\ \sin(\eta \mu(x_3)t) & \cos(\eta \mu(x_3)t) & 0 \\ 0 & 0 & 1
\end{matrix}
\right)
\end{equation}
with $\eta=5$ and $\mu(x_3) = \frac{1}{3}x_3^3$.
We append a right-hand side to  the equation (\ref{w10}) (and correspondingly to  (\ref{n10})) and set an initial value which are so that the exact solution of the transformed equation is in all cases $\tilde u(t,x)=e^{-6t}x_1x_2$. We calculated the parts from which the three right-hand sides (corresponding to the three cases) can be put together by hand and put them together within the implementation. For the discretization we use the coupling  $\Delta t=\Delta x^2$ between the step sizes in time and space.

We point out that the third example is chosen to observe the phenomenon of having only tangential motion which deteriorates the mesh in ESFEM (in the sense that the ratio of  the diameter and the incircle radius of a triangle might become large) and of course also affects the coefficients in the transformed equation.  We remark that for ESFEM the (relative) mesh size may change and tangential motion of the nodes may deteriorate the mesh. While the former might be compensated by appending additional nodes the latter can be compensated by introducing additional tangential motion of the mesh 
(ALE-ESFEM), see \cite{ES} for the latter. 

To be able to input the transformed coefficients as stated in (\ref{n10}) into the
Distributed
and Unified Numerics Environment (DUNE), see \cite{Bastian1, Bastian2, DUNE}, 
we present some auxiliary facts. Let $x=(x^i)_{1\le i\le 3}$ be Euclidean coordinates in $\mathbb{R}^3$ and $\xi=(\xi^j)_{1\le j\le 2}$ local coordinates in $\Gamma(0)$. Then $x=x^i(\xi^j)$ can be seen as a local representation of the embedding of $\Gamma(0)$ into $\mathbb{R}^3$ and $\Phi(x^i(\xi^j), t)$  of the embedding of $\Gamma(t)$ into $\mathbb{R}^3$. For the latter we drop the time dependence for simplicity and write 
$\Phi \circ x = \Phi(x^i(\xi^j))$.
We have
\begin{equation}
\begin{aligned}
g_{ij}(t) =& \frac{\partial}{\partial \xi^i}(\Phi \circ x) \cdot \frac{\partial}{\partial \xi^j}(\Phi \circ x) \\
\partial_kg_{ij}(t) =& \frac{\partial^2}{\partial \xi^i\partial \xi^k}(\Phi \circ x) \cdot \frac{\partial}{\partial \xi^j}(\Phi \circ x)
+ \frac{\partial^2}{\partial \xi^j\partial \xi^k}(\Phi \circ x) \cdot \frac{\partial}{\partial \xi^i}(\Phi \circ x),
\end{aligned}
\end{equation}
furthermore, 
\begin{equation}
\begin{aligned}
\nabla^{\Gamma(0)}_r g^{ij}(t) = -g^{ik}(t)g^{lj}(t)\partial_rg_{kl}(t)+ \Gamma(0)^i_{rl}g^{lj}+ \Gamma(0)^j_{rl}g^{li}.
\end{aligned}
\end{equation}
We calculate the tangential divergence of the evolution speed $v= \dot A(t)x$
\begin{equation}
\nabla^{\Gamma(t)} \cdot v = g^{ij}(t)(\frac{\partial}{\partial \xi^i}  \dot A(t)x+\dot A(t) \frac{\partial}{\partial \xi^i}x) \cdot \frac{\partial}{\partial \xi^j}(\Phi \circ x).
\end{equation}
Let $x=(x^i)\subset \Gamma(0)$ be given. Let $\alpha: \{1,2,3\}\rightarrow \{1,2,3\}$ be a bijection so that $|x^{\alpha(1)}| = \min_{j}|x^{\alpha(j)}|$, $x^{\alpha(3)}\neq 0$ and set
$x=x^{\alpha(3)}$, $y=x^{\alpha(2)}$, $z=x^{\alpha(1)}$,  $\varphi = \xi^2=\arctan(\frac{y}{x})$ and
\begin{equation}
\theta = \xi^1 = \begin{cases}\arccos(z) \text{ if } x  \ge 0\\
2 \pi - \arccos(z) \text{ else}.
\end{cases}
\end{equation}
Then we have
\begin{equation}
\begin{aligned}
x =& \sin \theta \cos \varphi \\
y =& \sin \theta \sin \varphi \\
z =& \cos \theta
\end{aligned}
\end{equation}
and
\begin{equation}
\begin{aligned}
\frac{\partial x}{\partial \xi ^1} =& \cos \theta \cos \varphi \\
\frac{\partial x}{\partial \xi ^2} =&  -\sin \theta \sin \varphi\\
\frac{\partial y}{\partial \xi ^1} =&  \cos \theta \sin \varphi \\
\frac{\partial y}{\partial \xi ^2} =&  \sin \theta \cos \varphi\\
\frac{\partial z}{\partial \xi ^1} =& -\sin \theta \\
\frac{\partial z}{\partial \xi ^2} =& 0. 
\end{aligned}
\end{equation}
We have to stipulate what we understand for a given triangle $T$ of $\Gamma_h(0)$ by a representation of $a^{ij}$ and $b^i$ in $\tilde T$ (where $\tilde T$ denotes the lift of $T$ to $\Gamma(0)$) with respect to given orthonormal basis $v_1, v_2$ in $T$. For this purpose we let $\tilde v_1, \tilde v_2$ be the orthogonal projections of $\frac{\partial}{\partial \xi^1}x$ and $\frac{\partial}{\partial \xi^2}x$ onto the plane which contains $T$. Then we define $\beta_k^ l$, $k,l=1,2$ by
\begin{equation}
\begin{aligned}
\tilde v_1 =& \beta_{1}^1 v_1 + \beta_{1}^2 v_2 \\
\tilde v_2 =& \beta_{2}^1 v_1 + \beta_{2}^2v_2
\end{aligned}
\end{equation} 
then
\begin{equation}
\tilde b^i = b^j \beta_{j}^i, \quad \tilde a^{ij} = a^{kl}\beta_{k}^i\beta_{l}^j
\end{equation}
are what we understand by evaluating $a^{ij}, b^i$ with respect to $v_1, v_2$. 
Let $P$ denote the plane containing $T$ and let $M$ be the matrix representation 
with respect to the standard basis in $\mathbb{R}^3$ of an (arbitrary) extension 
to $\mathbb{R}^3\times \mathbb{R}^3$
of the bilinear form on
$P\times P$ represented by the matrix $\tilde a^{ij}$ with respect to the basis
$v_1, v_2$. Correspondingly we construct a vector $B\in \mathbb{R}^3$ from $\tilde b^i$.
Then $M$ and $B$ are the coefficients which are compatible with the input format of DUNE.

In order to calculate the right-hand sides of the transformed equation for the exact solution $\tilde u$ (this means we evaluate the left-hand side of (\ref{n10}) for our special $\tilde u$) we use our coordinates $(\xi^i)$ to calculate the summand
\begin{equation}
g^{ij}(t)\nabla^{\Gamma(0)}_i\nabla^{\Gamma(0)}_j \tilde u  
\end{equation}
which can be put together from 
\begin{equation}
\nabla^{\Gamma(0)}_i\nabla^{\Gamma(0)}_j \tilde u = \frac{\partial^2 }{\partial \xi^i \partial \xi^j}\tilde u-\Gamma^k_{ij}(0)\frac{\partial }{\partial \xi^k} \tilde u.
\end{equation}
Note, that $\tilde u=e^{-6t}x_1x_2\neq e^{-6t}xy$ in general (depending on $\alpha$).

We use the mesh generator Gmsh, see \cite{Gmsh}, and for the implementation, more precisely,
DUNE-FEM, a discretization module for solving PDEs 
which depends on DUNE.

\end{document}